# Novel integral representations of the Riemann zeta-function and Dirichlet eta-function, close expressions for Laurent series expansions of powers of trigonometric functions and digamma function, and summation rules

Sergey K. Sekatskii

*Laboratoire de Physique de la Matière Vivante, IPHYS, BSP 408, Ecole Polytechnique Fédérale de Lausanne (EPFL), CH-1015 Lausanne, Switzerland.*

*E-mail: serguei.sekatski@epfl.ch*

We have established novel integral representations of the Riemann zeta-function and Dirichlet eta-function based on powers of trigonometric functions and digamma function, and then we use these representations to find close forms of Laurent series expansions of these same powers of trigonometric functions and digamma function. The so obtained series can be used to find numerous summation rules for certain values of the Riemann zeta and related functions and numbers, such as e.g. Bernoulli and Euler numbers.

## 1. Introduction

The aims of the present Note are the following. We start with establishing (an infinite number of) integral representations of the Riemann zeta-function and Dirichlet eta-function based on powers of trigonometric functions - for example, for

$0<\text{Re}s<2$, $\varsigma(s+2) = \dfrac{\pi^{s+1} \sin(\pi s/2)}{s+1} \int_0^\infty x^{-s-1}\left(\coth^2 x - \dfrac{1}{x^2} - \dfrac{2}{3}\right)dx$ and

$$\dfrac{1}{6}(s+2)(s+3)\varsigma(s+4) - \dfrac{4\pi^2}{3}\varsigma(s+2) =$$
$$\dfrac{\pi^{s+3} \sin(\pi s/2)}{s+1} \int_0^\infty x^{-s-1}\left(\coth^4 x - \dfrac{1}{x^4} - \dfrac{4}{3x^2} - \dfrac{26}{45}\right)dx.$$

(See e.g. [1] for the general discussion of the Riemann zeta-function and other special functions closely related to it). These representations continue a line of research which has been started in year 2000, and since then the following formulae were published respectively in [2 - 4]:

$$\varsigma(s) = \pi^{s-1} \sin(\pi s/2) \int_0^\infty \dfrac{1}{\tanh u}\left(\dfrac{1}{u^s} - \dfrac{1}{\sinh^s u}\right)du, \quad 1<\text{Re}s<2 \tag{1},$$

$$\varsigma(s+1) = \dfrac{\pi^s \cos(\pi s/2)}{2^{-s} - 1} \int_0^\infty x^{-s-1}\left(\dfrac{1}{\sinh x} - \dfrac{1}{x}\right)dx, \quad -1<\text{Re}s<1 \tag{2},$$

$$\varsigma(s) = -\dfrac{\pi^{s-1} \sin(\pi s/2)}{s-1} \int_0^\infty x^{1-s}\left(\dfrac{1}{\sinh^2 x} - \dfrac{1}{x^2}\right)dx, \quad 0<\text{Re}s<2 \tag{3}.$$

These representations, although mostly novel to the best of the author's knowledge (but see the Remark at the end of the Section 2.2), are however quite simple and almost trivial. Then we note that these same representations can be seen as Mellin transforms: for example, for (2) we have

$$\int_0^\infty x^{-s-1}\left(\dfrac{1}{\sinh x} - \dfrac{1}{x}\right)dx = \dfrac{(2^{-s} - 1)\pi^{-s}\varsigma(s+1)}{\cos(\pi s/2)} \tag{4},$$

and that the Mellin transform inversion theorem (see e.g. [5]) can be applied. Such an application readily gives close forms of Laurent series for powers of trigonometric functions, such as, for example,

$$\coth^4 x = \dfrac{1}{x^4} + \dfrac{4}{3x^2} + \dfrac{26}{45} +$$
$$\dfrac{2}{\pi^2}\sum_{n=1}^\infty (-1)^n (2n+1)(\dfrac{1}{6\pi^2}(2n+2)(2n+3)\varsigma(2n+4) - \dfrac{4}{3}\varsigma(2n+2))\left(\dfrac{x}{\pi}\right)^{2n} \tag{5}.$$

valid for $|x|<\pi$, and this subject seems less trivial. The multiplication of these series then leads to summation rules involving the values of the Riemann zeta-function, Dirichlet eta-function, or closely related with them Bernoulli or Euler numbers. The summation rules presented in the paper seem to be already known, but it might be that new rules can be obtained in this way.

In the last part of the paper we along the same lines briefly discuss integral representations based on powers of digamma (and polygamma) function as well as, again, corresponding close forms of Laurent series expansions for these functions, and following from their multiplication summation rules. These rules now include harmonic numbers.

To finish the introduction, let us also mention in passing that the point $s=0$ for (2), as well as the point $s=1$ for (3), are of no peculiarities. Simple analysis shows

$$\int_0^\infty x^{-1}\left(\frac{1}{\sinh x} - \frac{1}{x}\right)dx = -\ln 2 \text{ and } \int_0^\infty \left(\frac{1}{\sinh^2 x} - \frac{1}{x^2}\right)dx = -1.$$

## 2. Trigonometric functions
### 2.1. Proof of integral representation formulae

Let us consider a contour integral $\int_C z^{-s-1}\left(\cot z - \frac{1}{z}\right)dz$ taken round the contour $C$ consisting of the line $-iX$, $+iX$ and demi-circle connecting these points and lying in the right complex semi-plane. Here $X$ is real and tends to infinity, and in this limit for $0<\text{Re}\,s<1$ the contour integral is equal to ($z=ix$):

$$-\int_{-i\infty}^{i\infty} z^{-s-1}\left(\cot z - \frac{1}{z}\right)dz = -\int_{-\infty}^{\infty}(ix)^{-s-1}\left(\coth x - \frac{1}{x}\right)dx = -(i^{-s-1}-(-i)^{-s-1})\int_0^\infty x^{-s-1}\left(\coth x - \frac{1}{x}\right)dx =$$

$$2i\sin(\frac{s+1}{2}\pi)\int_0^\infty x^{-s-1}\left(\coth x - \frac{1}{x}\right)dx \tag{6}$$

(Initial minus sign comes from the necessity to round the contour counterclockwise).

From other point of view, this contour integral is equal, by Cauchy residue theorem, to absolutely converging sum $2\pi i \pi^{-s-1}\sum_{n=1}^\infty n^{-s-1} = 2i\pi^{-s}\varsigma(s+1)$. This finishes the proof of the representation

$$\varsigma(s+1) = \cos(\frac{s}{2}\pi)\pi^s \int_0^\infty x^{-s-1}\left(\coth x - \frac{1}{x}\right)dx \tag{7}$$

valid in the strip $0<\text{Re}\,s<1$. Absolutely in the same fashion $\int_C \tan z \cdot z^{-s-1}dz = 2i\sin(\frac{s+1}{2}\pi)\int_0^\infty \tanh x \cdot x^{-s-1}dx$, so that for the same strip

$\varsigma(s+1, 1/2) = \cos(\frac{s}{2}\pi)\pi^s \int_0^\infty \tanh x \cdot x^{-s-1}dx$ and

$$\varsigma(s+1) = \cos(\frac{s}{2}\pi)\frac{\pi^s}{2^{s+1}-1}\int_0^\infty \tanh x \cdot x^{-s-1}dx \tag{7a}$$

Relation $\varsigma(s, 1/2) = (2^s - 1)\varsigma(s)$ [1], where $\varsigma(s, a) = \sum_{n=0}^{\infty}(n+a)^{-s}$ is Hurwitz zeta-function, was used when deriving (7a).

Quite similarly, for the same strip we can consider a contour integral $\int_C z^{-s-1}\left(\frac{1}{\sin z} - \frac{1}{z}\right)dz$ and get $\cos(\pi s/2)\int_0^{\infty} x^{-s-1}\left(\frac{1}{\sinh x} - \frac{1}{x}\right)dx = -\pi^{-s}\eta(s+1)$, where $\eta(s, a) = \sum_{n=0}^{\infty}(-1)^{n-1}(n+a)^{-s}$ and we denote $\eta(s, 1) = \eta(s)$ - this is the Dirichlet's eta-function for which we know $\varsigma(s) = \frac{1}{1-2^{1-s}}\eta(s)$ [1]. Function $\eta(s, a) = \sum_{n=0}^{\infty}(-1)^{n-1}(n+a)^{-s}$ is a certain particular example of Lerch zeta-function $L(\lambda, a, s) = \sum_{n=0}^{\infty}\frac{\exp(2\pi i\lambda n)}{(n+a)^s}$; $L(1/2, a, s) = \eta(s, a)$. Thus we have proven the representation given in eq. (2) for this strip. This representation gives also an analytical continuation to the strip $-1 < \operatorname{Re} s \leq 0$.

Next we consider a contour integral $\int_C z^{-s-1}\left(\cot^2 z - \frac{1}{z^2} + \frac{2}{3}\right)dz$ in the strip 0<Re$s$<2. We have instead of (6)

$$-\int_{-i\infty}^{i\infty} z^{-s-1}\left(\cot^2 z - \frac{1}{z^2} + \frac{2}{3}\right)dz = -i\int_{-\infty}^{\infty}(ix)^{-s-1}\left(-\coth^2 x + \frac{1}{x^2} + \frac{2}{3}\right)dx =$$

$$-i(i^{-s-1} + (-i)^{-s-1})\int_0^{\infty} x^{-s-1}\left(\coth^2 x + \frac{1}{x^2} + \frac{2}{3}\right)dx = -2i\cos(\frac{s+1}{2}\pi)\int_0^{\infty} x^{-s-1}\left(-\coth^2 x + \frac{1}{x^2} + \frac{2}{3}\right)dx.$$

Application of Cauchy residue theorem (second order poles are located at $s=2n$, $n=1, 2, 3\ldots$) gives for the same integral $-2\pi i(\pi^{-s-2})\sum_{n=1}^{\infty}n^{-s-2}(s+1) = -2i\pi^{-s-1}(s+1)\varsigma(s+2)$, and thus we get

$$\varsigma(s+2) = \frac{\pi^{s+1}\sin(\pi s/2)}{s+1}\int_0^{\infty} x^{-s-1}\left(\coth^2 x - \frac{1}{x^2} - \frac{2}{3}\right)dx \qquad (8).$$

Using $\coth^2 x = \frac{1}{\sinh^2 x} + 1$, we can rewrite this as

$$\varsigma(s+2) = \frac{\pi^{s+1}\sin(\pi s/2)}{s+1}\int_0^{\infty} x^{-s-1}\left(\frac{1}{\sinh^2 x} - \frac{1}{x^2} + \frac{1}{3}\right)dx \qquad (9).$$

Quite similarly other trigonometric functions can be exploited, and here we will present only

$$\eta(s+1, 1/2) = -\pi^s \sin(\pi s/2)\int_0^{\infty} x^{-s-1}\left(\frac{1}{\cosh x} - 1\right)dx \qquad (10)$$

and
$$\varsigma(s+2) = \frac{\pi^{s+1} \sin(\pi s/2)}{(s+1)(2^{s+2}-1)} \int_0^\infty x^{-s-1} \left( \frac{1}{\cosh^2 x} - 1 \right) dx \qquad (11)$$
valid for 0<Re$s$<2.

Now, let us consider larger powers. For $\int_C z^{-s-1} \left( \cot^3 z - \frac{1}{z^3} + \frac{1}{z} \right) dz$ we first easily show that it is equal to $2i\cos(\pi s/2) \int_0^\infty x^{-s-1} \left( -\coth^3 x + \frac{1}{x^3} + \frac{1}{x} \right) dx$. Then we obtain for the same integral:

$$\int_C z^{-s-1} \left( \cot^3 z - \frac{1}{z^3} + \frac{1}{z} \right) dz = 2\pi i \pi^{-s-3} \left( \sum_{n=1}^\infty \frac{1}{2} n^{-s-3}(s+1)(s+2) - \pi^2 \sum_{n=1}^\infty n^{-s-1} \right)$$

$$= 2i\pi^{-s-2} \left( \frac{1}{2}(s+1)(s+2)\varsigma(s+3) - \pi^2 \varsigma(s+1) \right).$$ (For any $z = \pi n$, $n$=1, 2, 3... we, from Laurent series expansion $\cot^3 x = \frac{1}{x^3} - \frac{1}{x} + O(x)$ and periodicity of cotangent, have at this point a difference between a pole of the third order and a simple pole). Collecting everything together, we get for the strip 0<Re$s$<1:

$$\frac{1}{2}(s+1)(s+2)\varsigma(s+3) - \pi^2 \varsigma(s+1) = \pi^{s+2} \cos(\pi s/2) \int_0^\infty x^{-s-1} \left( \coth^3 x - \frac{1}{x^3} - \frac{1}{x} \right) dx \qquad (12).$$

*Stricto senso*, this cannot be named an integral representation of the Riemann zeta-function because two different species of such a function occur in the l.h.s. of (12). However, the relation of this result with the previous ones is evident and we will continue to use this terminology.

We have for $\int_C \tan^3 z \cdot z^{-s-1} dz$ that this is equal to $-2i\cos(\pi s/2) \int_0^\infty \tanh^3 x \cdot x^{-s-1} dx$, and further:

$$\int_C \tan^3 z \cdot z^{-s-1} dz = 2\pi i \pi^{-s-3} \left( \sum_{n=0}^\infty -\frac{1}{2}(n+1/2)^{-s-3}(s+1)(s+2) + \pi^2 \sum_{n=0}^\infty (n+1/2)^{-s-1} \right)$$

$$= 2i\pi^{-s-2} \left( -\frac{1}{2}(s+1)(s+2)\varsigma(s+3, 1/2) + \pi^2 \varsigma(s+1, 1/2) \right).$$ (For any $z = \pi(n+1/2)$, $n$=0, 1, 2, 3... we have at this point a difference between a simple pole and that of the third order). Now, using again $\varsigma(s, 1/2) = (2^s - 1)\varsigma(s)$, we have

$$-\frac{1}{2}(s+1)(s+2)(2^{s+3}-1)\varsigma(s+3) + \pi^2(2^{s+1}-1)\varsigma(s+1) = \pi^{s+2} \cos(\pi s/2) \int_0^\infty \tanh^3 x \cdot x^{-s-1} dx \quad (12a).$$

Quite similarly we, starting from the contour integral $\int_C z^{-s-1}\left(\cot^4 z - \dfrac{1}{z^4} + \dfrac{4}{3z^2} - \dfrac{26}{45}\right)dz$, find a representation

$$\dfrac{1}{6}(s+2)(s+3)\varsigma(s+4) - \dfrac{4\pi^2}{3}\varsigma(s+2) = \dfrac{\pi^{s+3}\sin(\pi s/2)}{s+1}\int_0^\infty x^{-s-1}\left(\coth^4 x - \dfrac{1}{x^4} - \dfrac{4}{3x^2} - \dfrac{26}{45}\right)dx$$

(13)

valid for 0<Re$s$<2.

Of course, similar procedure can be applied to numerous other functions like $\tan^n z$, $\sin^{-n} z$, $\cos^{-n} z$, $\cot^n z$, and others. For this, one needs to start with contour integrals $\int_C z^{-s-1}(\cot^n z - P_n(z))dz$ (or similar) where $P_n(z)$ is polynomial in non-positive powers of $z$ such that for $z$ tending to zero we respectively have $\cot^n z - P_n(z) = O(z)$ for odd $n$ and $\cot^n z - P_n(z) = O(z^2)$ for even $n$.

Note also that starting from contour integrals of the type $\int_C z^{-s-1}\left(\dfrac{1}{\sin(a+z)} + \dfrac{1}{\sin(a-z)} - \dfrac{2}{\sin a}\right)dz$, $\int_C z^{-s-1}\left(\dfrac{1}{\sin(a+z)} - \dfrac{1}{\sin(a-z)}\right)dz$ or similar we in the same manner as shown above obtain integral representation of the Hurwitz zeta function $\varsigma(s, a/\pi)$.

*2.2.    Laurent series expansions obtained from above integral representations*

Now we note that all appeared above integral representations are nothing else than Mellin transforms. For example, from (7) we have:

$$\int_0^\infty x^{-s-1}\left(\dfrac{1}{\tanh x} - \dfrac{1}{x}\right)dx = \dfrac{\pi^{-s}\varsigma(s+1)}{\cos(\pi s/2)} \qquad (14).$$

All conditions of the inversion Mellin theorem, see e.g. Theorem 12, p. 39 of Ref. [5], hold (integral absolutely converges in the strip 0<Re$s$<1 and function $\dfrac{1}{\tanh x} - \dfrac{1}{x}$ is continuous along the vertical line), thus we have

$$\dfrac{1}{\tanh x} - \dfrac{1}{x} = \dfrac{1}{2\pi i}\int_{c-i\infty}^{c+i\infty} x^s\left(\dfrac{\pi^{-s}\varsigma(s+1)}{\cos(s\pi/2)}\right)ds \qquad (15),$$

where 0<c<1. (Usually Mellin transform is defined as an integral including factor $x^{s-1}$ rather than $x^{-s-1}$ under the integral sign. Correspondingly, factor $x^{-s}$ rather than $x^s$ appears in the inverse Mellin transform. This is not essential and all necessary changes easily can be done).

For $|x|<\pi$, this integral, using the same contour C as above, can be calculated as a contour integral via Cauchy residue theorem as equal to

$$-2\pi i \frac{2}{\pi}\sum_{n=1}^{\infty}(-1)^{n}\varsigma(2n)\left(\frac{x}{\pi}\right)^{2n-1}$$ (the disappearance of the integral taken over demi-circle is evident; simple poles are located at $s=2n-1$, $n=1, 2, 3\ldots$; the initial minus comes from the necessity to round the contour counterclockwise). Thus we have:

$$\coth x = \frac{1}{x}+\frac{2}{\pi}\sum_{n=1}^{\infty}(-1)^{n-1}\left(\frac{x}{\pi}\right)^{2n-1}\varsigma(2n) \qquad (16),$$

which, of course, is the Laurent series expansion of the hyperbolic cotangent written usually in terms of Bernoulli numbers related with the Riemann zeta-function by $B_{2n}=\frac{(-1)^{n-1}2(2n)!}{(2\pi)^{2n}}\varsigma(2n)$ (see e.g. [6, 7] for basic properties of Bernoulli and Euler numbers; at the end of the paper we show how this relation can be obtained by our method): $\coth x = \frac{1}{x}+\sum_{n=1}^{\infty}\frac{2^{2n}B_{2n}}{(2n)!}x^{2n-1}=\frac{1}{x}+\frac{1}{3}x-\frac{1}{45}x^{3}+\frac{2}{945}x^{5}+\ldots$. From (7a) we have

$$\tanh x = \frac{2}{\pi}\sum_{n=1}^{\infty}(-1)^{n-1}(2^{2n}-1)\left(\frac{x}{\pi}\right)^{2n-1}\varsigma(2n) \qquad (17).$$

Completely similarly, starting from (3) rewritten as a Mellin transform, we get the relation

$$\frac{1}{\sinh x}=\frac{1}{x}+\frac{2}{\pi}\sum_{n=1}^{\infty}(-1)^{n}\eta(2n)\left(\frac{x}{\pi}\right)^{2n-1}=\frac{1}{x}+\frac{2}{\pi}\sum_{n=1}^{\infty}(-1)^{n}(1-2^{1-2n})\varsigma(2n)\left(\frac{x}{\pi}\right)^{2n-1} \qquad (18)$$

coinciding, of course, with the expansion

$$\frac{1}{\sinh x}=\frac{1}{x}-\sum_{n=1}^{\infty}\frac{2(2^{2n-1}-1)B_{2n}}{(2n)!}x^{2n-1}=\frac{1}{x}-\frac{1}{6}x+\frac{7}{360}x^{3}-\frac{31}{15120}x^{5}+\ldots \qquad (19)$$

and valid again for $|x|<\pi$. As well known, these series can be written as particular values of the hypergeometric function $\cot x = xF(1/2, 1; 3/2; -z^{2})$, $\frac{1}{\sin x}=xF(1/2, 1/2; 3/2; z^{2})$.

The same procedure applied to (8) rewritten as a Mellin transform gives $\coth^{2}x-\frac{1}{x^{2}}-\frac{2}{3}=\frac{1}{2\pi i}\int_{c-i\infty}^{c+i\infty}x^{s}\left(\frac{\pi^{-s-1}\varsigma(s+2)(s+1)}{\sin(s\pi/2)}\right)ds$ which, by Cauchy residue theorem leads to the Laurent series expansion

$$\coth^{2}x=\frac{1}{x^{2}}+\frac{2}{3}+\frac{2}{\pi^{2}}\sum_{n=1}^{\infty}(-1)^{n-1}(2n+1)\left(\frac{x}{\pi}\right)^{2n}\varsigma(2n+2) \qquad (20)$$

valid for $|x|<\pi$. Again using $\coth^{2}x=\frac{1}{\sinh^{2}x}+1$ we can rewrite as

$$\frac{1}{\sinh^{2}x}=\frac{1}{x^{2}}-\frac{1}{3}+\frac{2}{\pi^{2}}\sum_{n=1}^{\infty}(-1)^{n+1}(2n+1)\left(\frac{x}{\pi}\right)^{2n}\varsigma(2n+2) \qquad (21).$$

Similarly, we obtain from (10)

$$\frac{1}{\cosh x} = 1 + \frac{2}{\pi}\sum_{n=1}^{\infty}(-1)^n \left(\frac{x}{\pi}\right)^{2n} \eta(2n+1,\ 1/2) \tag{22}$$

and from (11)

$$\frac{1}{\cosh^2 x} = 1 + \frac{2}{\pi^2}\sum_{n=1}^{\infty}(-1)^n (2n+1)(2^{2n+2}-1)\left(\frac{x}{\pi}\right)^{2n} \varsigma(2n+2) \tag{23}.$$

These two expansions are valid for $|x| < \pi/2$. By definition of Euler numbers $E_{2n}$, Taylor expansion of hyperbolic secant is $\dfrac{1}{\cosh x} = \sum_{n=0}^{\infty}\dfrac{E_{2n}}{(2n)!}x^{2n} = 1 - \dfrac{x^2}{2} + \dfrac{5x^4}{24} - \dfrac{61x^6}{720} + ...,$ and we have an evident identification $E_{2n} = \dfrac{(-1)^n 2(2n)!}{\pi^{2n+1}} \eta(2n+1,\ 1/2)$.

Quite similarly, we get from (12) and (12a), the following Laurent series expansions valid for $|x| < \pi$ (24) or $|x| < \pi/2$ (25):

$$\coth^3 x = \frac{1}{x^3} + \frac{1}{x} + \frac{2}{\pi}\sum_{n=0}^{\infty}(-1)^{n+1}(\frac{1}{\pi^2}(n+1)(2n+3)\varsigma(2n+4) - \varsigma(2n+2))\left(\frac{x}{\pi}\right)^{2n+1} \tag{24}$$

$$\tanh^3 x = \frac{2}{\pi}\sum_{n=1}^{\infty}(-1)^{n+1}(\frac{1}{\pi^2}(n+1)(2n+3)(2^{2n+4}-1)\varsigma(2n+4) - (2^{2n+2}-1)\varsigma(2n+2))\left(\frac{x}{\pi}\right)^{2n+1} \tag{24a},$$

and from (13), for $|x| < \pi$

$$\coth^4 x = \frac{1}{x^4} + \frac{4}{3x^2} + \frac{26}{45} +$$

$$\frac{2}{\pi^2}\sum_{n=1}^{\infty}(-1)^n(2n+1)(\frac{1}{6\pi^2}(2n+2)(2n+3)\varsigma(2n+4) - \frac{4}{3}\varsigma(2n+2))\left(\frac{x}{\pi}\right)^{2n} \tag{25}.$$

This process of finding the Laurent series expansions can be continued indefinitely, as noted in the discussion at the end of the previous Section. Similarly to the note given in the end of the Section 2.1., we can establish Taylor expansions of the functions $\dfrac{1}{\sin^n(x+a)}$, $\dfrac{1}{\cos^n(x+a)}$ and like.

**Remark.** At this stage this is worthwhile to remind a number of somewhat similar "classic" representations of the Riemann zeta functions, such as $\int_0^{\infty}\dfrac{x^{\beta-1}}{\sinh(ax)}dx = \dfrac{2^{\beta}-1}{2^{\beta-1}a^{\beta}}\Gamma(\beta)\varsigma(\beta)$ (example N 3.523.1 of [8] book), $\int_0^{\infty}\dfrac{x^{\mu-1}e^{-\beta x}}{\sinh x}dx = 2^{1-\mu}\Gamma(\mu)\varsigma(\mu,\ \dfrac{1}{2}(\beta+1))$, N 3.552.1., etc. These Mellin transforms, as well as some others, are presented also in the classic Oberhettinger's "Tables of Melllin transforms" book [9], where we found, as an example 6.6, also $\int_0^{\infty}x^{z-1}(x^{-1} - \dfrac{1}{\sinh x})dx = 2(2^{-z}-1)\Gamma(z)\varsigma(z)$ valid for -1<Re$z$<1. (Note that this example is exactly eq. (2) presented above, which thus has been only reestablished in the 2007 paper [2]). To see the equivalence with our approach, one needs to use the functional equation [7] $\varsigma(s) = \frac{1}{\pi}(2\pi)^s \sin\frac{\pi s}{2}\Gamma(1-s)\varsigma(1-s)$ together with $\Gamma(1-s)\Gamma(s) = \frac{\pi}{\sin(\pi s)}$ [7] to obtain $\varsigma(s)\Gamma(s) = (2\pi)^s \varsigma(1-s)\frac{1}{2\cos(\pi s/2)}$, and remember that we consider integrals including factor $x^{-s-1}$ rather than $x^{s-1}$ under the integral sign.

### 2.3. Proof of summation formulae

Below we show how aforementioned Laurent series expansions can be used to establish numerous summation rules. We give rather short presentation because apparently all (or almost all) such summation rules seem are known, having being obtained by different powerful and often rather involved techniques, see e.g. [6 - 18] and references cited therein. (For example, an expression for the summation of quite general sums over products of Bernoulli numbers (but with the result expressed using inductively defined polynomials) has been presented by Huang and Huang [14]). Still we hope that an easy and transparent way to obtain such summation rules, and explicit form of the rules obtained can be useful. It appears also that among these summation rules including many combinations of the values of the Riemann zeta-functions as well as $\eta(s, 1/2)$ functions, etc. some still might be novel.

We know Laurent series expansion (20). On the other hand, we have by squaring eq. (16):

$$\coth^2 x = \frac{1}{x^2} + \frac{4}{\pi^2}\sum_{n=1}^{\infty}(-1)^{n-1}\left(\frac{x}{\pi}\right)^{2n-2}\varsigma(2n) + \frac{4}{\pi^2}\left(\sum_{n=1}^{\infty}(-1)^{n-1}\left(\frac{x}{\pi}\right)^{2n-1}\varsigma(2n)\right)^2 \qquad (26).$$

The comparison of the coefficients standing in front of the same powers of $(x/\pi)$ in (20), (26) gives the first summation rule: for $n=1, 2, 3\ldots$
$$2\sum_{k=1}^{n}\varsigma(2k)\varsigma(2n-2k+2) = (2n+3)\varsigma(2n+2).$$

For clarity of presentation, below we will formulate such rules as propositions.

**Proposition 1.** *For $n=1, 2, 3\ldots$, we have the summation rules*

$$\sum_{k=0}^{n}\varsigma(2k)\varsigma(2n-2k+2) = (n+1)\varsigma(2n+2) \qquad (27)$$

*and*

$$\sum_{k=0}^{n}C_{2n+2}^{2k}B_{2k}B_{2n-2k+2} = -(2n+2)B_{2n+2} \qquad (28).$$

Here and below we use binomial coefficients notation $C_n^k = \frac{n!}{k!(n-k)!}$. To get (27), we add also the term $2\varsigma(0)\varsigma(2n+2) = -\varsigma(2n+2)$ to both sides of the previous equality. Eq. (28) is just (27) written in the notation of Bernoulli numbers. This is, of course, Euler summation rule usually written as $\sum_{k=1}^{n} C_n^k B_k B_{n-k} = -nB_n$.

Next, squaring (18), comparing it with (21) and adding again the terms corresponding to $\varsigma(0)\varsigma(2n+2)$ to both sides of the equality, we arrive to the following proposition.

**Proposition 2.** *For n=1, 2, 3... we have the summation rules*

$$\sum_{k=0}^{n} (1-2^{1-2k})(1-2^{2k-2n-1})\varsigma(2k)\varsigma(2n-2k+2) = n\varsigma(2n+2) \qquad (29)$$

*and*

$$\sum_{k=0}^{n} C_{2n+2}^{2k}(1-2^{1-2k})(1-2^{2k-2n-1})B_{2n-2k+2}B_{2k} = -2nB_{2n+2} \qquad (30).$$

Eq. (30) is nothing else than the "nice sum identity"

$$\sum_{k=0}^{n} \frac{(1-2^{1-k})(1-2^{k-n+1})}{k!(n-k)!} B_{n-k}B_k = \frac{1-n}{n!} B_n \qquad (31)$$

attributed by Weisstein [6] to Gosper without further references (the present author failed to find an original proof of Gosper). To see this, one needs to add also a term corresponding to $B_0 B_{2n+2}$ to both sides of the above equation, and to note that odd Bernoulli numbers with $n>1$ are equal to zero while the $B_1=-1/2$ (or $+1/2$) term disappears due to the presence of the factors $1-2^{1-k}$, $1-2^{k+n-1}$. This identity is also nothing else than the relation $\sum_{k=0}^{2n} \frac{(2n)!}{k!(n-k)!} 2^k (-1)^{2n-k} B_{2n-k}B_k = (1-2n)B_{2n}$ obtained by Guinand using umbral method [17]. Moreover, its first appearance can be traced down to the famous Euler 1755 Memoire "Institutiones Calculi Differentialis".

Moving further, we multiply (16) and (20) and compare the result with (24). In such a way, we prove the following proposition.

**Proposition 3.** *For n=3,4,5... we have the summation rules*

$$4\sum_{m=1}^{n-2} (2m+1)\varsigma(2m+2)\varsigma(2n-2m-2) = 2(2n+1)(n-1)\varsigma(2n) - \frac{(2\pi)^2}{6}\varsigma(2n-2) \qquad (32)$$

*and*

$$\sum_{m=1}^{n-2} (2m+1)C_{2n}^{2m+2} B_{2m+2}B_{2n-2m-2} = -(2n+1)(n-1)B_{2n} + \frac{1}{24(2n-1)n} B_{2n-2} \qquad (33).$$

By squaring (20) and comparing with (25), we obtain:

**Proposition 4.** *For n=2,3,4... we have the summation rules*

$$\sum_{m=1}^{n-1}(2m+1)(2n-2m+1)\varsigma(2m+2)\varsigma(2n-2m+2) = \frac{1}{6}(2n+3)(n-1)(2n+5)\varsigma(2n+4) \quad (34)$$

*and*

$$\sum_{m=1}^{n-1}C_{2n+4}^{2m+2}(2m+1)(2n-2m+1)B_{2m+2}B_{2n-2m+2} = -\frac{1}{3}(2n+3)(n-1)(2n+5)B_{2n+4} \quad (35).$$

Again, this line of researches can be continued indefinitely long. Here, we finally illustrate our approach only with one more case establishing the simplest relation between Bernoulli and Euler numbers. (We do not count that relation between them which follows from the multiplication of Laurent series expansions of hyperbolic secant and cosecant and subsequent use of the formula $\frac{1}{\sinh x \cosh x} = \frac{2}{\sinh(2x)}$,

$$2\sum_{k=1}^{n}\eta(2k+1,\ 1/2)\eta(2n-2k+2) + \pi\eta(2n+3,\ 1/2) = \pi(2^{2n+1}-1)\eta(2n+2) \text{ for } n=1,\ 2,\ 3\ldots: \text{ the}$$

approach used in the present paper is superfluous to establish it). By squaring (22) and comparing it with (23), we prove

**Proposition 5.** *For n=1, 2, 3..., we have the summation rules*

$$2\sum_{k=1}^{n}\eta(2k+1,\ 1/2)\eta(2n-2k+3,\ 1/2) = (2n+3)(2^{2n+4}-1)\varsigma(2n+4) - 2\pi\eta(2n+3,\ 1/2) \quad (36)$$

*and*

$$\sum_{k=1}^{n}C_{2n+2}^{2k}E_{2k}E_{2n-2k+2} = \frac{2^{2n+4}(2^{2n+4}-1)}{2n+4}B_{2n+4} - 2E_{2n+2} \quad (37).$$

This is equivalent to the relation $(n+2)\sum_{k=0}^{n}\frac{n!}{k!(n-k)!}E_{n-k}E_k = 2^{n+2}(2^{n+2}-1)B_{n+2}$ obtained by Guinand using umbral method [17].

Other type of summation rules can be obtained using the products of the type $\sinh^m x/\cosh^n x$, $\cosh^m x/\sinh^n x$: they contain not the *products* of the Riemann zeta-functions (Bernoulli numbers) but *single such functions* (numbers). Here we illustrate this point with two examples following from our Taylor developments of the cubes of hyperbolic cotangent and tangent (24) and (24a). We have $\coth^3 x = \cosh^2 x \cdot \frac{\cosh x}{\sinh^3 x} = \cosh^2 x \cdot \frac{d}{dx}(-\frac{1}{2\sinh^2 x})$ and

$\tanh^3 x = \sinh^2 x \cdot \frac{\sinh x}{\cosh^3 x} = \sinh^2 x \cdot \frac{d}{dx}(-\frac{1}{2\cosh^2 x})$. Using the derivatives of our developments (21), (23), as well as the trivial $\cosh^2 x = \frac{1}{2} + \frac{1}{2}\sum_{n=0}^{\infty}\frac{(2x)^{2n}}{(2n)!}$ and $\sinh^2 x = \frac{1}{2}\sum_{n=1}^{\infty}\frac{(2x)^{2n}}{(2n)!}$, we easily obtain the following proposition, presented below only for Bernoulli numbers for the sake of brevity:

**Proposition 6.** *For n=2, 3, 4…, we have the summation rules*

$$\sum_{k=1}^{n-1} C_{2n+2}^{2k+2} k(2k+1)B_{2k+2} + 1 = (n+1)(2n+1)B_{2n} \qquad (38)$$

and

$$-\sum_{k=1}^{n-1} C_{2n+4}^{2k+2} k(2k+1)(2^{2k+2}-1)B_{2k+2} =$$

$$(2n+2)(2n+3)(2^{2n+4}-1)B_{2n+4} + (2^{2n+2}-1)(2n+3)(n+2)(2n^2+n+1)B_{2n+2} \qquad (39).$$

Finally, we would like to make the following remarks. The multiplication of Taylor developments of trigonometrical functions to establish summation rules dates back to at least 1824 when M. L. C. Bouvier used the comparison of the product of developments of *tan(x)* and *1/cos²(x)* with that of *tan(x)/cos²(x)* to establish the relation equivalent to our eq. (30) [18]. (Note that $d\tan x/dx = 1/\cos^2(x)$, $d(1/\cos^2(x))/dx = 2\tan x/\cos^2 x$ which enables to establish the necessary Taylor expansions by simple differentiation of the corresponding expansion of *tan(x)*; no any other work is needed. In this sense, the approach of the present paper looks essentially more general). Second, we need to note that the approach of Flajolet and Salvy paper [13] is especially close to that presented above. However, their approach cannot be used to directly obtain Laurent series expansions of the related functions.

### 3. Digamma function

Now let us briefly analyze what the same approach can bring to digamma function $\psi(z)$ defined as $\psi(z) = -\gamma + \sum_{n=0}^{\infty}\left(\frac{1}{n+1} - \frac{1}{n+z}\right)$; here $\gamma$ is Euler – Mascheroni constant; see e.g. [7] for main properties of the digamma function. (Indeed, in our setting we need to consider symmetrical or anti-symmetrical combinations of the type $\psi(x) \pm \psi(-x)$). Here we will consider only a few simplest examples.

The starting formulae are quite similar to the trigonometric function case. It is enough to know elementary properties of digamma function (easily following from definition) that in the vicinity of zero $\psi(z) = -\frac{1}{z} - \gamma + O(z)$, recurrence relation $\psi(z+1) = \psi(z) + \frac{1}{z}$, and the asymptotic $\psi(z) = O(\ln z)$ for large $z$, to start from the surely existing for Re$s>0$ contour integral $\int_C z^{-s-1}(\psi(-z) + \psi(z) + 2\gamma)dz$. For the strip 0<Re$s$<1, we immediately obtain a representation/Mellin transform

$$\varsigma(s+1) = -\frac{\pi}{\sin(\pi s/2)} \int_0^\infty x^{-s-1}(\psi(ix) + \psi(-ix) + 2\gamma)dx \qquad (40),$$

and inversion Mellin transform theorem (again, it is evident that all necessary conditions for this theorem hold) then readily gives Taylor expansion $\psi(ix) + \psi(-ix) + 2\gamma = 2\sum_{n=1}^\infty \varsigma(2n+1)x^{2n}$. Quite similarly, consideration of $\int_C z^{-s-1}(\psi(z) - \psi(-z) + \frac{2}{z})dz$ leads, for the same strip, to

$$i\cos(\pi s/2)\varsigma(s+1) = \pi \int_0^\infty x^{-s-1}\left(i\psi(ix) - i\psi(-ix) + \frac{2}{x}\right)dx \qquad \text{and then to}$$

$i\psi(ix) - i\psi(-ix) + \frac{2}{x} = 2\sum_{n=1}^\infty \varsigma(2n)x^{2n-1}$. Thus we have easily established (of course, well known) Laurent expansion

$$\psi(x) = -\frac{1}{x} - \gamma + \sum_{n=1}^\infty (-1)^{n-1}\varsigma(n+1)x^n \qquad (41)$$

valid for |x|<1. Its combination with the recurrence relation immediately supplies Laurent expansions in the vicinity of the simple poles $n=-1,-2,\ldots$

$$\psi(-n+x) = -\frac{1}{x} + H_n - \gamma + \sum_{k=1}^\infty (H_n^{(k+1)} + (-1)^{k+1}\varsigma(k+1))x^k \qquad (42)$$

where $H_n^{(k)} = \sum_{l=1}^n \frac{1}{l^k}$ is generalized harmonic number, $H_0^{(k)} = 0$ by definition. Correspondingly, in the vicinity of positive integers $n$, we have (by applying recurrence relation "in other direction")

$$\psi(n+x) = H_{n-1} - \gamma + \sum_{k=1}^\infty (-1)^k (H_{n-1}^{(k+1)} - \varsigma(k+1))x^k \qquad (43).$$

The series (42) converges for |x|<1 and that of eq. (43) for |x|<n.

The situation becomes much more involved already for the second power. Considering the simplest example, $\int_C z^{-s-1}\left(\psi(z)\psi(-z) + \frac{1}{z^2} - \gamma^2 - \frac{\pi^2}{3}\right)dz$, we do not get the integral representation of the Riemann zeta function but

$$\sin(\pi s/2)\int_0^\infty x^{-s-1}\left(\psi(ix)\psi(-ix) - \frac{1}{x^2} - \gamma^2 - \frac{\pi^2}{3}\right)dx = \pi\sum_{n=1}^\infty \psi(n)n^{-s-1} \qquad (44)$$

instead - due to the just established fact that at poles $x=1, 2, 3\ldots$ the residues are equal to $\psi(n) = H_{n-1} - \gamma$.

Let us *define* the function

$$\chi(s) = \sum_{n=1}^{\infty} \psi(n) n^{-s} = \sum_{n=1}^{\infty} H_{n-1} n^{-s} - \gamma\varsigma(s) = \sum_{n=1}^{\infty} H_n n^{-s} - \gamma\varsigma(s) - \varsigma(s+1).$$

Similar functions have been introduced and studied before, see, for example, [19] and references cited therein. This function enables us to write the corresponding Laurent expansion in compact form

$$\psi(x)\psi(-x) + \frac{1}{x^2} - \gamma^2 - \frac{\pi^2}{3} = 2\sum_{n=1}^{\infty} \chi(2n+1) x^{2n} \qquad (45).$$

Now multiplication of (39) taken for positive and negative $x$ and comparison with (43) give the first summation rule, which is the classical Euler summation rule, see e.g. Theorem 2.2 in [12]:

**Proposition 7.** *For odd $m= 3, 5, 7\ldots$, we have the summation rule*

$$\sum_{k=1}^{m-2} \varsigma(k+1)\varsigma(m-k) = (m+2)\varsigma(m+1) - 2\sum_{n=1}^{\infty} \frac{H_n}{n^m} \qquad (46).$$

Elaborating this line a bit further and comparing so obtained Taylor expansions for squares of $(\psi(x)+\psi(-x)+2\gamma)^2$, $(\psi(x)-\psi(-x)+2/x)^2$ with those directly following from the product of combinations of digamma function involved, we easily obtain again Euler summation rule (27) which does not contain any harmonic number, as well as the following proposition:

**Proposition 8.** *For odd $m= 2l+1=3, 5, 7\ldots$, we have the summation rule*

$$\sum_{k=1}^{l} \varsigma(2k+1)\varsigma(2l-2k+1) = 2l\varsigma(2l+2) - 2\sum_{n=1}^{\infty} \frac{H_n}{n^{2l+1}} \qquad (47).$$

*For even $m= 2l=4, 6, 8\ldots$, we have the summation rule*

$$\sum_{k=1}^{l} \varsigma(2k)\varsigma(2l-2k+1) = (2l+2)\varsigma(2l+1) - 2\sum_{n=1}^{\infty} \frac{H_n}{n^{2l}} \qquad (48).$$

Taken all together, they constitute classical Euler summation rule, see e.g. Theorem 2.2 in [12]:

**Proposition 9.** *For all $m=2, 3, 4\ldots$, we have the summation rule*

$$\sum_{k=1}^{m-2} \varsigma(k+1)\varsigma(m-k) = (m+2)\varsigma(m+1) - 2\sum_{n=1}^{\infty} \frac{H_n}{n^m} \qquad (49).$$

As for the case of trigonometric functions, this line of researches can be continued much further. In particular, polygamma functions and "shifted" functions of the type $\psi(x+n)$ can be useful.

Finally, to add certain "self-consistence" to the paper, let us underline that the relation between Bernoulli numbers and some particular values of the Riemann zeta-functions, $B_{2n} = \frac{(-1)^{n-1} 2(2n)!}{(2\pi)^{2n}} \varsigma(2n)$, can be easily established by the presented approach if one will define Bernoulli numbers as $\frac{t}{e^t - 1} = \sum_{n=0}^{\infty} \frac{B_n}{n!} t^n$ and start from the contour integrals $\int_C z^{-s-1} \left( \frac{iz}{e^{iz} - 1} + \frac{iz}{e^{-iz} - 1} \right) dz$ and $\int_C z^{-s-1} \left( \frac{iz}{e^{iz} - 1} - \frac{iz}{e^{-iz} - 1} - 2 \right) dz$. This consideration is completely similar to the handling of the functions $\psi(x) \pm \psi(-x)$ presented above in this Section.

**Acknowledgement**

The author thanks Jacques Gélinas for useful remarks, suggestions and additional references.